\documentclass[10pt, a4paper] {article}

\usepackage{amssymb}
\usepackage{amsthm}
\usepackage{amsmath}
\usepackage{amsfonts}
\usepackage{mathtools}
\usepackage[shortlabels]{enumitem}
\usepackage[sharp]{easylist}

\theoremstyle{definition}
\newtheorem{definition}{Definition}

\theoremstyle{proposition}
\newtheorem{proposition}{Proposition}

\theoremstyle{lemma}
\newtheorem{lemma}{Lemma}

\theoremstyle{corollary}
\newtheorem{corollary}{Corollary}

\theoremstyle{affirmation}
\newtheorem{affirmation}{Affirmation}

\newenvironment{prova}{\setlength{\parindent}{0pt}{\textbf{Proof}.}}{\hspace{\stretch{1}}$\blacksquare$}

\title{PARTIAL MODEL THEORY - ULTRAPRODUCTS AND COMPACTNESS}

\author{Rodolfo Cunha Carnier\\
\footnotesize University of São Paulo (BR)}

\date{}

\begin{document}

\maketitle

\begin{abstract}
In the present paper we prove the compactness theorem with respect to partial structures and quasi-truth, using the technique of ultraproducts. Partial structures and quasi-truth are two notions developed within the partial structures approach, which is a philosophical conception that emerged in the context of contemporary philosophy of science. Nevertheless, the notions developed within this conception, in particular the two mentioned, have a model-theoretic content that has not been explored so far, so that this paper is part of a project where we intend to analyze their formal properties by means of the development of a partial model theory, which is an extension of traditional model theory to partial structures.
\end{abstract}

                                  \section{Introduction}

The present paper is part of a project where we develop what we call \textit{partial model theory}, which is an extension of traditional model theory to \textit{partial structures}. These are a specific kind of structure, developed within the partial structures approach, which is a philosophical conception that constitutes the \textit{semantic approach of theories}. One of the main aims of the partial structures approach is to supply a conceptual frame that allows the use of logical concepts, in contexts where our knowledge about a domain under investigation is incomplete. And in addition to the notion of partial structure, there are two other main notions developed within the partial structures approach that sustain it and contribute to its purposes; they are the notions of \textit{partial relation} and \textit{quasi-truth}.

In summary, an \textit{n}-ary partial relation defined over a set \textit{A} is a relation that is not defined for every \textit{n}-tuple of objects of \textit{A}, whereas a partial structure is a structure whose relations defined over its universe are partial. The idea is that a partial structure models a determinate domain of investigation from an epistemic point of view, based on what is known and unknown about it. Thus if we are analyzing a domain of investigation and we do not know whether or not certain entities stand in a relation to each other, we can formally accommodate this situation by employing a partial structure. It is to be remarked that partial structures play the same role as that of usual structures, in the sense that they also provide an interpretation for the symbols of a language, so we can talk about \textit{partial models for a language} (cf. \cite{Bueno1996} p. 191).

As to quasi truth, from a philosophical point of view it may be considered a \textit{pragmatic} truth notion, since it was conceived with the purpose of apprehending some of the ``intentions'' underlying the truth notions of Charles Sanders Peirce and William James (cf. \cite{DaCosta2003} p. 12); accordingly, quasi-truth can also be called \textit{pragmatic truth}. From a logical point of view, however, quasi-truth was inspired by the standard notion of truth of Alfred Tarski, which is reflected in some formal characteristics that both notions share, like the fact that a sentence can only be quasi-true in a partial structure with respect to an interpretation, the same way a sentence can only be true (in the Tarskian sense) in a structure with respect to an interpretation.

The partial structures approach emerged in the context of contemporary philosophy of science, and its concepts have been used notably in this field (cf. \cite{DaCosta1990} and \cite{Bueno1997})\footnote{In \cite{Mikenberg1986}, nonetheless, the conceptual frame of the partial structures approach is used to solve a mathematical problem.}. But those concepts have a model-theoretic content that has not been explored so far and, as we have argued elsewhere\footnote{We are mentioning a paper to appear.}, their use in the philosophy of science and in other areas as well could be improved by means of the development of that model-theoretic content. In our aforementioned work we put forward such a content introducing the basics of partial model theory, so that in the present paper we intend to give a step further, by showing that this theory preserves an important result of traditional model theory; more precisely, we will prove the \textit{compactness theorem} with respect to partial structures and quasi-truth, using the technique of \textit{ultraproducts}.

                                  \section{Partial model theory: a minimum}

This section might be seen as a prolegomenon to the main part of this paper, contained in the next section. In what follows we will introduce the most basic notions and results of partial model theory, that we shall use to prove the compactness theorem. Some of these notions will be introduced in a particular way, so as to overcome certain criticisms addressed to the partial structures approach. Therefore, by developing our partial model theory we are also reformulating the partial structures approach with the aim of circumventing the referred criticisms.

\subsection{Partial relations, partial functions and partial structures}

              \begin{definition}
Let \textit{A} be a set. An \textit{n}-ary \textit{partial relation} \textbf{R} over \textit{A} is an ordered triple (\textbf{R}\textsubscript{+}, \textbf{R}\textsubscript{-}, \textbf{R}\textsubscript{0}) such that:

\begin{enumerate}
\item \textbf{R}\textsubscript{+}, \textbf{R}\textsubscript{-}, \textbf{R}\textsubscript{0} are mutually disjoint sets;

\item $\textbf{R}\textsubscript{+} \cup \textbf{R}\textsubscript{-} \cup \textbf{R}\textsubscript{0} = A^n$.
\end{enumerate}

In case $\textbf{R}\textsubscript{0} = \emptyset$, \textbf{R} is a usual \textit{n}-ary relation that may be identified with \textbf{R}\textsubscript{+}, and is said to be a \textit{total relation}.
              \end{definition}

The idea is that \textbf{R} is a relation which is not defined for every \textit{n}-tuple of objects of \textit{A}, so that \textbf{R}\textsubscript{+} is the set of the \textit{n}-tuples that satisfy \textbf{R}, \textbf{R}\textsubscript{-} is the set of the \textit{n}-tuples that do not satisfy \textbf{R}, and \textbf{R}\textsubscript{0} is the set of the \textit{n}-tuples for which it is not defined whether or not they satisfy \textbf{R}.

              \begin{definition}
Let \textit{A} be a set. An \textit{n}-ary \textit{partial function} \textit{f} over \textit{A} is an $(n+1)$-ary partial relation $(f\textsubscript{+}, f\textsubscript{-}, f\textsubscript{0})$ over \textit{A}, such that for every $a_1, ..., a_n, b \in A$:
\begin{enumerate}
\item If $(a_1, ..., a_n, b) \in f\textsubscript{+}$, then $(a_1, ..., a_n, b') \in f\textsubscript{-}$ for every $b' \in A$ such that $b \neq b'$;
\item If $(a_1, ..., a_n, b) \in f\textsubscript{-}$, then there exists $b' \in A$ such that either $(a_1, ..., a_n, b') \in f\textsubscript{+}$ or $(a_1, ..., a_n, b') \in f\textsubscript{0}$;
\item If $(a_1, ..., a_n, b) \in f\textsubscript{0}$, then for every $b' \in A$ either $(a_1, ..., a_n, b') \in f\textsubscript{-}$ or $(a_1, ..., a_n, b') \in f\textsubscript{0}$.
\end{enumerate}

In case $f\textsubscript{0} = \emptyset$, \textit{f} is a usual \textit{n}-ary function and is said to be a \textit{total function}.
              \end{definition}

As can be seen, if \textit{f} is an \textit{n}-ary partial function over a set \textit{A} according to Definition 2, then \textit{f} is a partial relation $(f\textsubscript{+}, f\textsubscript{-}, f\textsubscript{0})$. Moreover, if $a_1, ..., a_n, b \in A$ are any elements and $(a_1, ..., a_n, b) \in f\textsubscript{+}$, then it is defined that $(a_1, ..., a_n, b)$ satisfies \textit{f}, and for every $b' \in A$ it is defined that $(a_1, ..., a_n, b')$ does not satisfy \textit{f}. If $(a_1, ..., a_n, b) \in f\textsubscript{-}$ then it is defined that $(a_1, ..., a_n, b)$ does not satisfy \textit{f}, and there exists $b' \in B$ such that either it is defined that $(a_1, ..., a_n, b')$ satisfies \textit{f} or it is not defined whether $(a_1, ..., a_n, b')$ satisfies \textit{f}. Finally, if $(a_1, ..., a_n, b) \in f\textsubscript{0}$ then it is not defined whether $(a_1, ..., a_n, b)$ satisfies \textit{f}, and for each $b' \in B$, either it is not defined whether $(a_1, ..., a_n, b')$ satisfies \textit{f} or it is defined that $(a_1, ..., a_n, b')$ does not satisfy \textit{f}.

This notion of partial function was developed to account for a criticism against the notion of partial structure, made by Sebastian Lutz. According to this criticism, partial structures cannot accommodate the sort of lack of knowledge about functions typically found in science, because they assign function symbols to partial functions in the usual mathematical sense (cf. \cite{Lutz2021} pp. 1357-8). The problem with a partial structure assigning function symbols to this kind of partial function is that in some contexts we may not know the value of a function for a given argument, but we may know that certain elements of the function's counterdomain do not correspond to that value. Nevertheless, if \textit{f} is an \textit{n}-ary partial function defined over a set \textit{A} and the value of \textit{f} for an argument $(a_1, ..., a_n)$ is not defined, then we cannot say whether or not $f(a_1, ..., a_n) \neq a$, for each $a \in A$. Thus, in order to formally accommodate the sort of lack of knowledge about functions that Lutz has in mind, one needs a partial function such that even if its value is not defined for a certain argument, there are certain elements of its couterdomain that cannot correspond to the value in question. It can be easily seen that our new notion of partial function behaves precisely that way.

              \begin{definition}
A \textit{language} is a pair $\mathcal{L} = (L, (\mu, \delta))$ such that \textit{L} is the union of two disjoint sets $\mathcal{F}$ and $\mathcal{R}$, while $(\mu, \delta)$ is a pair of mappings $\mu : \mathcal{F} \to \mathbb{N}$ and $\delta : \mathcal{R} \to \mathbb{N}^*$. The elements of $\mathcal{F}$ and $\mathcal{R}$ are called \textit{function symbols} and \textit{relation symbols}, respectively. We also assume that $\mathcal{L}$ is composed by a set of the usual \textit{logical symbols} such as \textit{variables}, \textit{operators}, \textit{quantifiers}, the \textit{identity sign}, and \textit{brackets}. Further, we assume that $\mathcal{L}$ is composed by a set $\mathcal{C} \subseteq \mathcal{F}$, such that $\mu(\mathcal{C}) = \{0\}$ and whose members are called \textit{constant symbols}. If $\mathcal{R}$ is empty then $\mathcal{L}$ is said an \textit{algebraic language}, and if $\mathcal{F}$ is empty then $\mathcal{L}$ is said a \textit{relational language}. The pair $(\mu, \delta)$ is the \textit{type} of $\mathcal{L}$.
              \end{definition}

              \begin{definition}
Let $\mathcal{L} = (L, (\mu, \delta))$ be a language. An $\mathcal{L}$-\textit{partial structure} $\mathcal{A}$ is a pair $\mathcal{A} = (A, (Z^\mathcal{A})_{Z \in L})$ where \textit{A} is a non-empty set and the family $(Z^\mathcal{A})_{Z \in L}$ is such that:

\begin{enumerate}
\item If $Z \in \mathcal{C}$, then $Z^\mathcal{A} \in A$ whenever $Z^\mathcal{A}$ is defined;
\item If $Z \in \mathcal{F}$, then $Z^\mathcal{A}$ is a $\mu(Z)$-ary partial function over \textit{A};
\item If $Z \in \mathcal{R}$, then $Z^\mathcal{A}$ is a $\delta(Z)$-ary partial relation over \textit{A}.
\end{enumerate}
              \end{definition}

The set \textit{A} is the \textit{universe} of the structure $\mathcal{A}$ and $Z^\mathcal{A}$ is the \textit{interpretation} of the symbol \textit{Z} in $\mathcal{A}$. Notice that if \textit{Z} is a constant symbol, $Z^\mathcal{A}$ is a 0-ary partial function over $\mathcal{A}$, so that if $Z^\mathcal{A}$ is total then $Z^\mathcal{A} \in A$, whilst if $Z^\mathcal{A}$ is not total we let \textit{Z} be a constant symbol whose interpretation is not defined\footnote{Our main motivation for defining partial structures like this, that is, in such a way that they may not interpret constants, is another criticism of Sebastian Lutz. According to this criticism, partial structures cannot accommodate the sort of lack of knowledge about constants typically found in science, because they interpret every constant symbol of a language, so that in some contexts we may not know the value of a given constant, at least up to a certain instant of time (cf. \cite{Lutz2021} pp. 1356-7). It is easy to see, nonetheless, that our new notion of partial structure is not subject to such a criticism - nor to the criticism regarding partial functions. As a conclusion, we have that by employing this new notion of partial structure together with our new notion of partial function, lack of knowledge about constants and functions can be easily accommodate within the conceptual frame of the partial structures approach.}. If for each $Z \in \mathcal{F} \cup \mathcal{R}$ it follows that $Z^\mathcal{A}$ is total, $\mathcal{A}$ is said to be a \textit{total} $\mathcal{L}$-\textit{structure} (or simply \textit{total structure}) - note that total structures are usual structures with which we work in traditional model theory. When $\mathcal{L}$ is an algebraic language $\mathcal{A}$ is said a \textit{partial algebra}, and when $\mathcal{L}$ is a relational language $\mathcal{A}$ is said a \textit{relational partial structure}. The cardinal of $\mathcal{A}$ is the cardinal of its universe, i.e., $Card(\mathcal{A}) = |A|$.

Two $\mathcal{L}$-partial structures $\mathcal{A} = (A, (Z^\mathcal{A})_{Z \in L})$ and $\mathcal{B} = (B, (Z^\mathcal{B})_{Z \in L})$ are equal if:
\begin{enumerate} [(a)]
\item $A = B$;
\item For $Z \in \mathcal{C}$, we have:
\begin{enumerate} [i.]
\item $Z^\mathcal{A}$ is defined if and only if $Z^\mathcal{B}$ is defined;
\item If $Z^\mathcal{A}$ and $Z^\mathcal{B}$ are both defined, $Z^\mathcal{A} = Z^\mathcal{B}$.
\end{enumerate}

\item For $Z \in \mathcal{F} \cup \mathcal{R}$, we have that $Z^\mathcal{A} = Z^\mathcal{B}$.
\end{enumerate}

To conclude, notice that every partial structure might be converted into a relational partial structure - just like every structure can be converted into a relational structure in traditional model theory (cf. \cite{Manzano1999} p. 7). Thus we will not lose generality by restricting ourselves to relational partial structures. So, from now on, we will assume a fixed relational language $\mathcal{L}$, and we shall work with the class of (relational) $\mathcal{L}$-partial structures.

                                           \subsection{Expansion}

Together with the notions of partial structure and quasi-truth, the next notion that shall be introduced is characteristic of partial model theory. 

              \begin{definition}
Let $\mathcal{A} = (A, (Z^\mathcal{A})_{Z \in L})$ and $\mathcal{B} = (B, (Z^\mathcal{B})_{Z \in L})$ be $\mathcal{L}$-partial structures. We say that $\mathcal{B}$ \textit{expands} $\mathcal{A}$, in symbols $\mathcal{A} \Subset \mathcal{B}$, if:

\begin{enumerate}
\item $A = B$;
\item For $Z \in \mathcal{R}$, we have:
\begin{enumerate}
\item $Z^\mathcal{A}\textsubscript{+} \subseteq Z^\mathcal{B}$\textsubscript{+};
\item $Z^\mathcal{A}\textsubscript{-} \subseteq Z^\mathcal{B}$\textsubscript{-}.
\end{enumerate}
\end{enumerate}

              \end{definition}

Intuitively, the notion of expansion represents a possible increase of knowledge about the domain modeled by $\mathcal{A}$, so that $\mathcal{B}$ models the same domain taking this increase of knowledge into account.

\begin{proposition}
For each $\mathcal{L}$-partial structure $\mathcal{A}$ there exists an $\mathcal{L}$-partial structure expanding it.
\end{proposition}

\begin{prova}
In what follows, we will show that given an $\mathcal{L}$-partial structure $\mathcal{A} = (A, (Z^\mathcal{A})_{Z \in L})$, there are three $\mathcal{L}$-partial structures, which are possibly the same, that expand $\mathcal{A}$. The first of them is $\mathcal{A}$ itself, given that:

\begin{enumerate} [I.]
\item $A = A$;

\item For $Z \in \mathcal{R}$, it is immediate that $Z^\mathcal{A}\textsubscript{+} \subseteq Z^\mathcal{A}\textsubscript{+}$ and $Z^\mathcal{A}\textsubscript{-} \subseteq Z^\mathcal{A}\textsubscript{-}$. 
\end{enumerate} 

Hence, $\mathcal{A} \Subset \mathcal{A}$. Now define the pairs $\mathcal{B} = (B, (Z^\mathcal{B})_{Z \in L})$ and $\mathcal{D} = (D, (Z^\mathcal{D})_{Z \in L})$ as follows:

\begin{enumerate} [(1)]
\item $B = D = A$;

\item For $Z \in \mathcal{R}$, we have: 

\begin{enumerate} [i.]
\item $Z^\mathcal{B}\textsubscript{+} = Z^\mathcal{A}\textsubscript{+} \cup Z^\mathcal{A}\textsubscript{0}$ and $Z^\mathcal{D}\textsubscript{-} = Z^\mathcal{A}\textsubscript{-} \cup Z^\mathcal{A}\textsubscript{0}$;
\item $Z^\mathcal{B}\textsubscript{-} = Z^\mathcal{A}\textsubscript{-}$ and $Z^\mathcal{D}\textsubscript{+} = Z^\mathcal{A}\textsubscript{+}$;
\item $Z^\mathcal{B}\textsubscript{0} = Z^\mathcal{D}\textsubscript{0} = \emptyset$.
\end{enumerate}
\end{enumerate}

By construction, both $\mathcal{B}$ and $\mathcal{D}$ are total $\mathcal{L}$-structures. Furthermore, by construction we also have that $\mathcal{A} \Subset \mathcal{B}$ and $\mathcal{A} \Subset \mathcal{D}$.
\end{prova}

\begin{corollary}
For each $\mathcal{L}$-partial structure $\mathcal{A}$:
\begin{enumerate}
\item $\mathcal{A} \Subset \mathcal{A}$;
\item There exists a total structure $\mathcal{B}$ such that $\mathcal{A} \Subset \mathcal{B}$.
\end{enumerate}
\end{corollary}

\begin{prova}
Straightforward.
\end{prova}

\begin{lemma}
Let $\mathcal{A} = (A, (Z^\mathcal{A})_{Z \in L})$ and $\mathcal{B} = (B, (Z^\mathcal{B})_{Z \in L})$ be $\mathcal{L}$-partial structures such that $\mathcal{A} \Subset \mathcal{B}$. For $Z \in \mathcal{R}$, it follows that if $Z^\mathcal{A}$ is total then $Z^\mathcal{A} = Z^\mathcal{B}$.
\end{lemma}

\begin{prova}
Assume that $Z \in \mathcal{R}$ is such that $Z^\mathcal{A}$ is total and let us check that $Z^\mathcal{A} = Z^\mathcal{B}$. We already have that $Z^\mathcal{A}\textsubscript{0} = \emptyset \subseteq Z^\mathcal{B}\textsubscript{0}$. Further, since $\mathcal{A} \Subset \mathcal{B}$, we also have that $Z^\mathcal{A}\textsubscript{+} \subseteq Z^\mathcal{B}\textsubscript{+}$ and $Z^\mathcal{A}\textsubscript{-} \subseteq Z^\mathcal{B}\textsubscript{-}$. Therefore, we just have to verify that $Z^\mathcal{B}\textsubscript{+} \subseteq Z^\mathcal{A}\textsubscript{+}$, $Z^\mathcal{B}\textsubscript{-} \subseteq Z^\mathcal{A}\textsubscript{-}$ and $Z^\mathcal{B}\textsubscript{0} \subseteq Z^\mathcal{A}\textsubscript{0}$. Thus, let $a_1, ..., a_{\delta(Z)} \in A$ be such that $(a_1, ..., a_{\delta(Z)}) \in Z^\mathcal{B}\textsubscript{+}$. Given that $Z^\mathcal{A}\textsubscript{0} = \emptyset$, $Z^\mathcal{A}\textsubscript{-} \subseteq Z^\mathcal{B}\textsubscript{-}$ and $Z^\mathcal{B}\textsubscript{+} \cap Z^\mathcal{B}\textsubscript{-} = \emptyset$, clearly $(a_1, ..., a_{\delta(Z)}) \in Z^\mathcal{A}\textsubscript{+}$, so that $Z^\mathcal{A}\textsubscript{+} = Z^\mathcal{B}\textsubscript{+}$. Using an analogous argument it follows that $Z^\mathcal{A}\textsubscript{-} = Z^\mathcal{B}\textsubscript{-}$. But then, it is immediate that $Z^\mathcal{B}\textsubscript{0} = Z^\mathcal{A}\textsubscript{0}$ and hence $Z^\mathcal{A} = Z^\mathcal{B}$.
\end{prova}

             \begin{definition}
Let $\mathcal{A} = (A, (Z^\mathcal{A})_{Z \in L})$ and $\mathcal{B} = (B, (Z^\mathcal{B})_{Z \in L})$ be $\mathcal{L}$-partial structures. We say that $\mathcal{B}$ is an $\mathcal{A}$-\textit{normal} structure (or simply $\mathcal{A}$-\textit{normal}), if $\mathcal{A} \Subset \mathcal{B}$ and $\mathcal{B}$ is total.
              \end{definition}

\begin{proposition}
Let $\mathcal{A}$ be an $\mathcal{L}$-partial structure. If $\mathcal{A}$ is a total structure then $\mathcal{A}$ is $\mathcal{A}$-normal, and for every $\mathcal{L}$-partial structure $\mathcal{B}$, if $\mathcal{B}$ is $\mathcal{A}$-normal then $\mathcal{A} = \mathcal{B}$.
\end{proposition}

\begin{prova}
Immediate by item 1 of Corollary 1 and by Lemma 1.
\end{prova}

\begin{proposition}
For each $\mathcal{L}$-partial structure $\mathcal{A}$ there exists an $\mathcal{A}$-normal structure. 
\end{proposition}

\begin{prova}
Immediate by item 2 of Corollary 1.
\end{prova}

\bigskip

The proof of Proposition 1 shows that for every $\mathcal{L}$-partial structure $\mathcal{A} = (A, (Z^{\mathcal{A}})_{Z \in L})$, there exist two distinctive $\mathcal{A}$-normal structures $\mathcal{B} = (B, (Z^\mathcal{B})_{Z \in L})$ and $\mathcal{D} = (D, (Z^\mathcal{D})_{Z \in L})$, such that for every relation symbol \textit{Z}, if it is not defined whether a given $\delta(Z)$-tuple $(a_1, ..., a_{\delta(Z)})$ satisfies $Z^\mathcal{A}$ then it is defined that $(a_1, ..., a_{\delta(Z)})$ satisfies $Z^\mathcal{B}$ and does not satisfy $Z^\mathcal{D}$. Those structures will be important in the subsequent discussion, so it will be convenient to have a specific notation for them. Accordingly, given an $\mathcal{L}$-partial structure $\mathcal{A} = (A, (Z^{\mathcal{A}})_{Z \in L})$, we shall use the symbols $\mathcal{A}_+$ and $\mathcal{A}_-$ to denote the $\mathcal{A}$-normal structures $\mathcal{B} = (B, (Z^\mathcal{B})_{Z \in L})$ and $\mathcal{D} = (D, (Z^\mathcal{D})_{Z \in L})$ respectively, such that for $Z \in \mathcal{R}$:

\begin{enumerate} [i.]
\item $Z^\mathcal{B}\textsubscript{+} = Z^\mathcal{A}\textsubscript{+} \cup Z^\mathcal{A}\textsubscript{0}$ and $Z^\mathcal{D}\textsubscript{-} = Z^\mathcal{A}\textsubscript{-} \cup \ Z^\mathcal{A}\textsubscript{0}$;
\item $Z^\mathcal{B}\textsubscript{-} = Z^\mathcal{A}\textsubscript{-}$ and $Z^\mathcal{D}\textsubscript{+} = Z^\mathcal{A}\textsubscript{+}$.
\end{enumerate}
 
\noindent Moreover, we will use the notation $A_+$ to denote the universe of $\mathcal{A}_+$ and $A_-$ to denote the universe of $\mathcal{A}_-$.

\bigskip

\begin{lemma}
Let $\mathcal{A} = (A, (Z^\mathcal{A})_{Z \in L})$ be an $\mathcal{L}$-partial structure and $\mathcal{B} = (B, (Z^\mathcal{B})_{Z \in L})$ be an $\mathcal{A}$-normal structure. For $Z \in \mathcal{R}$, it follows that $Z^\mathcal{B}\textsubscript{+} \subseteq Z^\mathcal{A}\textsubscript{+} \cup Z^\mathcal{A}\textsubscript{0}$ and $Z^\mathcal{B}\textsubscript{-} \subseteq Z^\mathcal{A}\textsubscript{-} \cup Z^\mathcal{A}\textsubscript{0}$.
\end{lemma}

\begin{prova}
Let $a_1, ..., a_{\delta(Z)} \in A$ be such that $(a_1, ..., a_{\delta(Z)}) \in Z^\mathcal{B}\textsubscript{+}$. So $(a_1, ..., a_{\delta(Z)}) \not\in Z^\mathcal{B}\textsubscript{-}$, and since $Z^\mathcal{A}\textsubscript{-} \subseteq Z^\mathcal{B}\textsubscript{-}$ it follows that $(a_1, ..., a_{\delta(Z)}) \not\in Z^\mathcal{A}\textsubscript{-}$. But then $(a_1, ..., a_{\delta(Z)}) \in A^{\delta(Z)} - Z^\mathcal{A}\textsubscript{-}$ and hence $(a_1, ..., a_{\delta(Z)}) \in Z^\mathcal{A}\textsubscript{+} \cup Z^\mathcal{A}\textsubscript{0}$. The argument to show that $Z^\mathcal{B}\textsubscript{-} \subseteq Z^\mathcal{A}\textsubscript{-} \cup Z^\mathcal{A}\textsubscript{0}$ is similar.
\end{prova}

\subsection{Quasi-truth}

We conclude this section introducing the notion of quasi-truth and some other related notions, like the notions of \textit{quasi-validity}, \textit{quasi-consequence} and \textit{quasi-equivalence}. We will also make some remarks about the logic underlying quasi-truth.

              \begin{definition}
Let $\alpha$ be an $\mathcal{L}$-sentence and $\mathcal{A}$ an $\mathcal{L}$-partial structure. We say that $\alpha$ is \textit{quasi-true in} $\mathcal{A}$, in symbols $\mathcal{A} \mid\models \alpha$, if there exists an $\mathcal{A}$-normal structure $\mathcal{B}$ such that $\mathcal{B} \models \alpha$, i.e., if $\alpha$ is true in $\mathcal{B}$ in the Tarskian sense. Otherwise, $\alpha$ is said to be \textit{quasi-false in} $\mathcal{A}$.
              \end{definition}

              \begin{definition}
An $\mathcal{L}$-sentence $\alpha$ is said \textit{quasi-valid}, in symbols $\mid\models \alpha$, if for every $\mathcal{L}$-partial structure $\mathcal{A}$ it follows that $\mathcal{A} \mid\models \alpha$.
              \end{definition}

              \begin{definition}
If $\mathcal{A} \mid\models \alpha$, then $\mathcal{A}$ is said to be a \textit{partial model} of $\alpha$. Given a set $\Gamma$ of $\mathcal{L}$-sentences, we say that $\mathcal{A}$ is a \textit{partial model} of $\Gamma$, in symbols $\mathcal{A} \mid\models \Gamma$, if for every $\gamma \in \Gamma$ we have that $\mathcal{A} \mid\models \gamma$.
              \end{definition}

              \begin{definition}
An $\mathcal{L}$-sentence $\alpha$ is a \textit{(logical)} \textit{quasi-consequence} of an $\mathcal{L}$-sentence $\gamma$, in symbols $\gamma \mid\models \alpha$, if every partial model of $\gamma$ is a partial model of $\alpha$. An $\mathcal{L}$-sentence $\alpha$ is a \textit{(logical)} \textit{quasi-consequence} of a set of $\mathcal{L}$-sentences $\Gamma$, in symbols $\Gamma \mid\models \alpha$, if every partial model of $\Gamma$ is a partial model of $\alpha$.
              \end{definition}

              \begin{definition}
Two $\mathcal{L}$-sentences $\alpha$ and $\beta$ are \textit{(logically)} \textit{quasi-equivalent} if $\alpha \mid\models \beta$ and $\beta \mid\models \alpha$.
              \end{definition}

The next result shows that quasi-truth is equivalent to Tarskian truth when restricted to total structures.

\begin{proposition}
Let $\mathcal{A}$ be a total $\mathcal{L}$-structure and $\varphi$ an $\mathcal{L}$-sentence. Then $\mathcal{A} \mid\models \varphi$ if and only if $\mathcal{A} \models \varphi$.
\end{proposition}

\begin{prova}
$(\Rightarrow)$ Suppose that $\mathcal{A} \mid\models \varphi$. Then, there exists an $\mathcal{A}$-normal structure $\mathcal{B}$ such that $\mathcal{B} \models \varphi$. But $\mathcal{A} = \mathcal{B}$ according to Proposition 2. Hence, $\mathcal{A} \models \varphi$.

\smallskip

$(\Leftarrow)$ Now suppose that $\mathcal{A} \models \varphi$. Using again Proposition 2 we have that $\mathcal{A}$ is $\mathcal{A}$-normal and therefore  
$\mathcal{A} \mid\models \varphi$.
\end{prova}

Using Proposition 4, the reader may verify that the notions of quasi-validity and quasi-equivalence coincide with the notions of validity and equivalence respectively. The next example shows, nonetheless, that the notion of quasi-consequence does not coincide with the notion of consequence.

Suppose that $\mathcal{L} = (L, (\mu, \delta))$ is such that: $L = \mathcal{R} = \{R\}$ and $\delta(R) = 1$. Now assume that $\mathcal{A} = (A, (Z^\mathcal{A})_{Z \in L})$ is an $\mathcal{L}$-partial structure where: $A = \{a\}$, $R^\mathcal{A}\textsubscript{0} = \{a\}$ and $R^\mathcal{A}\textsubscript{+} = R^\mathcal{A}\textsubscript{-} = \emptyset$. So, the only $\mathcal{A}$-normal structures are the structures $\mathcal{A}_+$ and $\mathcal{A}_-$, which are such that: $R^{\mathcal{A}_+}\textsubscript{+} = R^{\mathcal{A}_-}\textsubscript{-} = \{a\}$ and $R^{\mathcal{A}_+}\textsubscript{-} = R^{\mathcal{A}_+}\textsubscript{0} = R^{\mathcal{A}_-}\textsubscript{+} = R^{\mathcal{A}_-}\textsubscript{0} = \emptyset$. But then, clearly $\mathcal{A}_+ \models \forall x Rx$ and $\mathcal{A}_- \models \neg\forall x Rx$, so that $\mathcal{A} \mid\models \forall x Rx$ and $\mathcal{A} \mid\models \neg\forall x Rx$. On the other hand, $\mathcal{A} \mid\not\models \forall x Rx \wedge \neg\forall x Rx$, because $\mathcal{A}_+ \not\models \forall x Rx \wedge \neg\forall x Rx$ and $\mathcal{A}_- \not\models \forall x Rx \wedge \neg\forall x Rx$. Therefore, assuming that $\varphi$ is $\forall x Rx$ and $\psi$ is $\forall x Rx \wedge \neg\forall x Rx$, it follows that

\begin{center}
$\{\varphi, \neg\varphi\} \mid\not\models \psi$.
\end{center}

This example also implies the next result which will be stated as an affirmation.

\begin{affirmation}
The following assertions are true:

\begin{enumerate}
\item If $\varphi$ is $\forall x Rx$ and $\psi$ is $\neg\forall x Rx$, then $\{\varphi, \psi\} \mid\not\models \varphi \wedge \psi$;

\item If $\varphi$ is $\forall x Rx$ and $\psi$ is $\forall x Rx \wedge \neg\forall x Rx$, then $\{\varphi \to \psi, \varphi\} \mid\not\models \psi$ and $\{\varphi \lor \psi, \neg\varphi\} \mid\not\models \psi$;

\item If $\varphi$ is $\neg(\forall x Rx \wedge \neg\forall x Rx)$ and $\psi$ is $\forall x Rx$, then $\{\varphi \to \psi, \neg\psi\} \mid\not\models \neg\varphi$;

\item If $\varphi$ is $\neg\forall x Rx$, $\psi$ is $\forall x Rx \wedge \neg\forall x Rx$ and $\Gamma$ is \{$\forall x Rx$\}, then $\Gamma \mid\models \varphi \to \psi$ and $\Gamma \cup \{\varphi\} \mid\not\models \psi$.
\end{enumerate}

\end{affirmation}

Therefore the notion of quasi-consequence does not coincide with the notion of consequence, and thus the logic underlying quasi-truth is not classical logic.

                                  \section{Ultraproducts and the compactness theorem}

In order to prove the compactness theorem via ultraproducts, it is necessary to introduce some preliminary notions and results. So, in the next pages we shall proceed as follows: we will begin by defining the concept of \textit{direct product} of a family of partial structures and by demonstrating some propositions related to the concepts of total relation, total structure and quasi-truth. Next we will define the concepts of \textit{filter} and \textit{reduced product}, as well as the concepts of \textit{ultrafilter} and ultraproduct, ending with a series of results of which the compactness theorem will be the last. As we shall see, our proof of the latter is identical to that of traditional model theory, so that the most interesting features of what will be presented lie in what comes before the compactness theorem, notwithstanding the importance of the fact that it holds as regards partial structures and quasi-truth.

\subsection{Direct products} 

              \begin{definition}
Let \textit{I} be a non-empty set and \{$\mathcal{A}_i: i \in I$\} be a family of $\mathcal{L}$-partial structures. The \textit{direct product} (or simply \textit{product}) of \{$\mathcal{A}_i: i \in I$\} is an $\mathcal{L}$-partial structure $\mathcal{A}_I$ described as follows:

\begin{enumerate}

\item The universe of $\mathcal{A}_I$ is

\begin{center}
$A_I := \prod_{i \in I}A_i$ = \{$h : h$ is a mapping with $Dom(h) = I$ and $h(i) \in A_i$\};
\end{center}

\item For $Z \in \mathcal{R}$ and $h_1, ..., h_{\delta(Z)} \in A_I$, we have:

\begin{enumerate} [i.]
\item $(h_1, ..., h_{\delta(Z)}) \in Z^{\mathcal{A}_I}\textsubscript{+}$ if and only if $\{i \in I : (h_1(i), ..., h_{\delta(Z)}(i)) \in Z^{\mathcal{A}_i}\textsubscript{+}\} = I$;

\item $(h_1, ..., h_{\delta(Z)}) \in Z^{\mathcal{A}_I}\textsubscript{-}$ if and only if $\{i \in I : (h_1(i), ..., h_{\delta(Z)}(i)) \in Z^{\mathcal{A}_i}\textsubscript{-}\} \neq \emptyset$;

\item $(h_1, ..., h_{\delta(Z)}) \in Z^{\mathcal{A}_I}\textsubscript{0}$ if and only if $\{i \in I : (h_1(i), ..., h_{\delta(Z)}(i)) \in Z^{\mathcal{A}_i}\textsubscript{0}\} \neq \emptyset = \{i \in I : (h_1(i), ..., h_{\delta(Z)}(i)) \in Z^{\mathcal{A}_i}\textsubscript{-}\}$.
\end{enumerate}
\end{enumerate}

The structure $\mathcal{A}_i$ is called the \textit{i}th \textit{factor} of the product. If all partial structures $\mathcal{A}_i$ are equal to some fixed $\mathcal{L}$-partial structure $\mathcal{A}$, the product is said a \textit{direct power} (or simply \textit{power}) of $\mathcal{A}$, and we write $\mathcal{A}^I$. 
              \end{definition}

The idea behind the last item of clause 12.2 is that for $Z \in \mathcal{R}$, a $\delta(Z)$-tuple $(h_1, ..., h_{\delta(Z)})$ belongs to $Z^{\mathcal{A}_I}\textsubscript{0}$ if and only if for some $i \in I$, it is not defined whether $(h_1(i), ..., h_{\delta(Z)}(i))$ belongs to $Z^{\mathcal{A}_i}$, and for each $k \neq i$ either it is not defined whether $(h_1(k), ..., h_{\delta(Z)}(k))$ belongs to $Z^{\mathcal{A}_k}$ or it is defined that $(h_1(k), ..., h_{\delta(Z)}(k))$ indeed belongs to $Z^{\mathcal{A}_k}$.

Before proceeding, we make three further remarks. First, recall that whenever $Z^\mathcal{A}\textsubscript{0} = \emptyset$, given an $\mathcal{L}$-partial structure $\mathcal{A}$ and a relation symbol $Z \in \mathcal{R}$, the partial relation $Z^\mathcal{A}$ can be identified with $Z^\mathcal{A}\textsubscript{+}$ (see page 2 above); hence when $ Z^\mathcal{A}\textsubscript{+} = Z^\mathcal{A}\textsubscript{0} = \emptyset$, $Z^\mathcal{A}$ can be identified with the empty relation $\emptyset$. Second, if $\mathcal{A} = (A, (Z^\mathcal{A})_{Z \in L})$ and $\mathcal{B} = (B, (Z^\mathcal{B})_{Z \in L})$ are total $\mathcal{L}$-structures and we want to prove that $\mathcal{A} = \mathcal{B}$, then for $Z \in \mathcal{R}$ we just have to verify that $Z^\mathcal{A}\textsubscript{+} = Z^\mathcal{B}\textsubscript{+}$, because since $Z^\mathcal{A}\textsubscript{0} = Z^\mathcal{B}\textsubscript{0} = \emptyset$, in case $Z^\mathcal{A}\textsubscript{+} = Z^\mathcal{B}\textsubscript{+}$ we will also have that $Z^\mathcal{A}\textsubscript{-} = Z^\mathcal{B}\textsubscript{-}$ and therefore $Z^\mathcal{A} = Z^\mathcal{B}$. From now on, we will use this fact implicitly. Finally, in what follows, unless otherwise stated we shall work with a fixed direct product $\mathcal{A}_I$ of a fixed family \{$\mathcal{A}_i: i \in I$\} of $\mathcal{L}$-partial structures, where \textit{I} is non-empty.

\begin{proposition}
Consider the product $\mathcal{A}_I$ of the family \{$\mathcal{A}_i: i \in I$\}. For $Z \in \mathcal{R}$, it follows that $Z^{\mathcal{A}_I}\textsubscript{0} \neq \emptyset$ if and only if for some $i \in I$, $Z^{\mathcal{A}_i}\textsubscript{0} \neq \emptyset$ and for each $k \in I$, $Z^{\mathcal{A}_k} \neq \emptyset$.
\end{proposition}

\begin{prova}
($\Rightarrow$) Immediate by clause 12.2 of Definition 12.

\smallskip

($\Leftarrow$) Suppose that $Z^{\mathcal{A}_i}\textsubscript{0} \neq \emptyset$ for some $i \in I$, and for each $k \in I$ we have $Z^{\mathcal{A}_k} \neq \emptyset$. Thus for every \textit{k}, there exists a $\delta(Z)$-tuple $(a_{1k}, ..., a_{\delta(Z)k})$ such that $(a_{1k}, ..., a_{\delta(Z)k}) \in Z^{\mathcal{A}_k}\textsubscript{+} \cup \ Z^{\mathcal{A}_k}\textsubscript{0}$. In case $k = i$, there exists a $\delta(Z)$-tuple $(a_{1i}, ..., a_{\delta(Z)i}) \in Z^{\mathcal{A}_i}\textsubscript{0}$. Now if for every \textit{k} we select one of those tuples, we see that for each \textit{n}-coordinate $(1 \leq n \leq \delta(Z))$ there is a mapping, say $h_n \in A_I$, such that $h_n(k) = a_{nk}$. But then the mappings $h_1, ..., h_{\delta(Z)} \in A_I$ are such that for every $k \neq i$, $(h_1(k), ..., h_{\delta(Z)}(k)) \in Z^{\mathcal{A}_k}\textsubscript{+} \cup \ Z^{\mathcal{A}_k}\textsubscript{0}$, and for $k = i$, $(h_1(i), ..., h_{\delta(Z)}(i)) \in Z^{\mathcal{A}_i}\textsubscript{0}$. So, $(h_1, ..., h_{\delta(Z)}) \in Z^{\mathcal{A}_I}\textsubscript{0}$ and hence $Z^{\mathcal{A}_I}\textsubscript{0} \neq \emptyset$.
\end{prova}

\begin{corollary}
$\mathcal{A}_I$ is total if and only if for every $i \in I$, $\mathcal{A}_i$ is total, or whenever $Z^{\mathcal{A}_i}\textsubscript{0} \neq \emptyset$ there is $k \in I$ such that $Z^{\mathcal{A}_k} = \emptyset$.
\end{corollary}

\begin{prova}
Immediate by Proposition 5.
\end{prova}

These last two results entail that if $\mathcal{A}_I$ is total, then there are three possibilities: either (i) $\mathcal{A}_i$ is total for every $i \in I$, or (ii) there are $i, k \in I$ such that $\mathcal{A}_i$ is total and $\mathcal{A}_k$ is not, or yet, (iii) $\mathcal{A}_i$ is not total for each $i \in I$. On the other hand, if $\mathcal{A}_I$ is not total, then there is at least one partial structure $\mathcal{A}_i$ which is not total. But it may also happen that for some $i, k \in I$, $\mathcal{A}_i$ is total and $\mathcal{A}_k$ is not, or $\mathcal{A}_i$ is not total for each $i \in I$.

\begin{proposition}
Consider the family \{$\mathcal{A}_i: i \in I$\} and let the family \{$\mathcal{B}_i: i \in I$\} be such that $\mathcal{B}_i$ is $\mathcal{A}_i$-normal for each $i \in I$. Then $\mathcal{B}_I$ is $\mathcal{A}_I$-normal.
\end{proposition}

\begin{prova}
It is immediate that $\mathcal{B}_I$ is total, so we just have to check that $\mathcal{A}_I \Subset \mathcal{B}_I$. Thus:

\begin{enumerate} [I.]
\item Let us show that $A_I = B_I$.
\begin{eqnarray*}
h \in A_I &\Leftrightarrow& Dom(h) = I \ \text{and} \ h(i) \in A_i \hspace{0,9 cm} \text{(for each} \ i \in I)\\
&\Leftrightarrow& Dom(h) = I \ \text{and} \ h(i) \in B_i \hspace{0,89 cm} \text{(for each} \ i \in I)\\
&\Leftrightarrow& h \in B_I
\end{eqnarray*}

Hence, $A_I = B_I$.

\item For $Z \in \mathcal{R}$ and $h_1, ..., h_{\delta(Z)} \in A_I$:

\begin{enumerate} [ i.]
\item Assume that $(h_1, ..., h_{\delta(Z)}) \in Z^{\mathcal{A}_I}\textsubscript{+}$. So $\{i \in I : (h_1(i), ..., h_{\delta(Z)}(i)) \in Z^{\mathcal{A}_i}\textsubscript{+}\} = I$. Next let us verify that $\{i \in I : (h_1(i), ..., h_{\delta(Z)}(i)) \in Z^{\mathcal{B}_i}\textsubscript{+}\} = I$. It is immediate that $\{i \in I : (h_1(i), ..., h_{\delta(Z)}(i)) \in Z^{\mathcal{B}_i}\textsubscript{+}\} \subseteq I$. Now let $k \in I$ be any element. Then we have that $k \in \{i \in I : (h_1(i), ..., h_{\delta(Z)}(i)) \in Z^{\mathcal{A}_i}\textsubscript{+}\}$, so that $(h_1(k), ..., h_{\delta(Z)}(k)) \in Z^{\mathcal{A}_k}\textsubscript{+} \subseteq Z^{\mathcal{B}_k}\textsubscript{+}$. Hence $k \in \{i \in I : (h_1(i), ..., h_{\delta(Z)}(i)) \in Z^{\mathcal{B}_i}\textsubscript{+}\}$ and thus $\{i \in I : (h_1(i), ..., h_{\delta(Z)}(i)) \in Z^{\mathcal{B}_i}\textsubscript{+}\} = I$. Therefore, $(h_1, ..., h_{\delta(Z)}) \in Z^{\mathcal{B}_I}\textsubscript{+}$.

\item Suppose that $(h_1, ..., h_{\delta(Z)}) \in Z^{\mathcal{A}_I}\textsubscript{-}$. So $\{i \in I : (h_1(i), ..., h_{\delta(Z)}(i)) \in Z^{\mathcal{A}_i}\textsubscript{-}\} \neq \emptyset$, whence there exists $k \in I$ such that $(h_1(k), ..., h_{\delta(Z)}(k)) \in Z^{\mathcal{A}_k}\textsubscript{-}$. But then $(h_1(k), ..., h_{\delta(Z)}(k)) \in Z^{\mathcal{B}_k}\textsubscript{-}$ and thus $\{i \in I : (h_1(i), ..., h_{\delta(Z)}(i)) \in Z^{\mathcal{B}_i}\textsubscript{-}\} \neq \emptyset$. Hence, $(h_1, ..., h_{\delta(Z)}) \in Z^{\mathcal{B}_I}\textsubscript{-}$.
\end{enumerate}
\end{enumerate}

Therefore, $\mathcal{A}_I \Subset \mathcal{B}_I$ and so $\mathcal{B}_I$ is $\mathcal{A}_I$-normal.
\end{prova}

\begin{proposition}
Consider the product $\mathcal{A}_I$ of the family \{$\mathcal{A}_i: i \in I$\}. Then, ${\mathcal{A}_I}_+$ is the direct product of the family \{${\mathcal{A}_i}_+ : i \in I$\} and ${\mathcal{A}_I}_-$ is the direct product of the family \{${\mathcal{A}_i}_- : i \in I$\}.
\end{proposition}

\begin{prova}
Let $\mathcal{B} = (B, (Z^\mathcal{B})_{Z \in L})$ be the direct product of \{${\mathcal{A}_i}_+ : i \in I$\}, $\mathcal{D} = (D, (Z^\mathcal{D})_{Z \in L})$ be the direct product of \{${\mathcal{A}_i}_- : i \in I$\}, and let us check that ${\mathcal{A}_I}_+ = \mathcal{B}$ and ${\mathcal{A}_I}_- = \mathcal{D}$. 
We begin by the first equality. 

\begin{enumerate} [(a)]

\item The argument to show that ${A_I}_+ = B$ is similar to that of item I of Proposition 6.

\item For $Z \in \mathcal{R}$ and $h_1, ..., h_{\delta(Z)} \in {A_I}_+$, we have:
\begin{eqnarray*}
(h_1, ..., h_{\delta(Z)}) \in Z^{{\mathcal{A}_I}_+}\textsubscript{+} &\Leftrightarrow& (h_1, ..., h_{\delta(Z)}) \in Z^{\mathcal{A}_I}\textsubscript{+} \cup Z^{\mathcal{A}_I}\textsubscript{0}\\
&\Leftrightarrow& (h_1(i), ..., h_{\delta(Z)}(i)) \in Z^{\mathcal{A}_i}\textsubscript{+} \cup Z^{\mathcal{A}_i}\textsubscript{0} \hspace{1,2 cm} \text{(for each} \ i \in I)\\
&\Leftrightarrow& (h_1(i), ..., h_{\delta(Z)}(i)) \in Z^{{\mathcal{A}_i}_+}\textsubscript{+} \hspace{2,18 cm} \text{(for each} \ i \in I)\\
&\Leftrightarrow& (h_1, ..., h_{\delta(Z)}) \in Z^{\mathcal{B}}\textsubscript{+}
\end{eqnarray*}

Hence, $Z^{{\mathcal{A}_I}_+}\textsubscript{+} = Z^\mathcal{B}\textsubscript{+}$.
\end{enumerate}

Therefore, ${\mathcal{A}_I}_+ = \mathcal{B}$. Using an analogous argument it follows that ${\mathcal{A}_I}_- = \mathcal{D}$.
\end{prova}

\subsection{Reduced products}

In this part of section 3 we introduce the notion of reduced product, but first let us define the notion of filter.

\begin{definition}
Consider the set \textit{I}. A \textit{filter F over I} is defined to be a set of subsets of \textit{I} such that:

\begin{enumerate}
\item $I \in F$ and $\emptyset \not\in F$;

\item If $X, Y \in F$, then $X \cap Y \in F$;

\item If $X \in F$ and $X \subseteq Y \subseteq I$, then $Y \in F$.
\end{enumerate}
\end{definition}

\begin{affirmation}
Consider the set \textit{I}. Then the set $\{I\}$ is a filter over \textit{I}.  
\end{affirmation}

The filter $\{I\}$ is called the \textit{trivial filter} (cf. \cite{Chang1990} p. 211), and it will be very important for us later.

Consider the family \{$\mathcal{A}_i: i \in I$\} and let \textit{F} be a collection of subsets of \textit{I}. We define the relation $\sim_F$ over $A_I$ by

\begin{center}
$u \sim_F v$ if and only if $\{i \in I : u(i) = v(i)\} \in F$.
\end{center}

\begin{affirmation}
If \textit{F} is a filter over \textit{I}, then $\sim_F$ is an equivalence relation on $A_I$. Moreover, if \textit{F} is the trivial filter, $\sim_F$ is the identity relation on $A_I$.
\end{affirmation}

We write $[u]_F$ for the equivalence class of the element $u \in A_I$. Thus

\begin{center}
$[u]_F = \{v \in A_I : u \sim_F v\}$.
\end{center}

\noindent We also use the symbol $A_I/F$ to denote the set of all equivalence classes of $\sim_F$, that is,

\begin{center}
$A_I/F = \{[u]_F : u \in A_I\}$.
\end{center}

\begin{definition}
Let \textit{F} be a filter over \textit{I} and consider the relation $\sim_F$. We define the \textit{reduced product} of the family \{$\mathcal{A}_i: i \in I$\} over \textit{F}, in symbols $\mathcal{A}_I/F$, as follows: 

\begin{enumerate}
\item The universe of $\mathcal{A}_I/F$ is $A_I/F$;

\item For $Z \in \mathcal{R}$ and $[h_1]_F, ..., [h_{\delta(Z)}]_F \in A_I/F$, we have: 
\begin{enumerate} [i.]
\item $([h_1]_F, ..., [h_{\delta(Z)}]_F) \in Z^{\mathcal{A}_I/F}\textsubscript{+}$ if and only if $\{i \in I : (h_1(i), ..., h_{\delta(Z)}(i)) \in Z^{\mathcal{A}_i}\textsubscript{+}\} \in F$;

\item $([h_1]_F, ..., [h_{\delta(Z)}]_F) \in Z^{\mathcal{A}_I/F}\textsubscript{0}$ if and only if $\{i \in I : (h_1(i), ..., h_{\delta(Z)}(i)) \in Z^{\mathcal{A}_i}\textsubscript{0}\} \in F$;

\item $([h_1]_F, ..., [h_{\delta(Z)}]_F) \in Z^{\mathcal{A}_I/F}\textsubscript{-}$ otherwise.
\end{enumerate}
\end{enumerate}

If for each $i \in I$, $\mathcal{A}_i = \mathcal{A}$, where $\mathcal{A}$ is any $\mathcal{L}$-partial structure, the reduced product is denoted by $\mathcal{A}^I/F$ and is said a \textit{reduced power} of $\mathcal{A}$ over \textit{F}.
\end{definition}

We next check that clause 2 of Definition 14 does not depend on the representatives $h_1, ..., h_{\delta(Z)}$, but only on the equivalence classes $[h_1]_F, ..., [h_{\delta(Z)}]_F$.

\begin{proposition}
Let \textit{F} be a filter over \textit{I}, and consider the relation $\sim_F$ as well as the reduced product $\mathcal{A}_I/F$. If $Z \in \mathcal{R}$, $[h_1]_F, ..., [h_{\delta(Z)}]_F, [h'_1]_F, ..., [h'_{\delta(Z)}]_F \in A_I/F$ and $h_1 \sim_F h'_1, ..., h_{\delta(Z)} \sim_F h'_{\delta(Z)}$, then:

\begin{enumerate}
\item $([h_1]_F, ..., [h_{\delta(Z)}]_F) \in Z^{\mathcal{A}_I/F}\textsubscript{+}$ if and only if $([h'_1]_F, ..., [h'_{\delta(Z)}]_F) \in Z^{\mathcal{A}_I/F}\textsubscript{+}$;

\item $([h_1]_F, ..., [h_{\delta(Z)}]_F) \in Z^{\mathcal{A}_I/F}\textsubscript{0}$ if and only if $([h'_1]_F, ..., [h'_{\delta(Z)}]_F) \in Z^{\mathcal{A}_I/F}\textsubscript{0}$;

\item $([h_1]_F, ..., [h_{\delta(Z)}]_F) \in Z^{\mathcal{A}_I/F}\textsubscript{-}$ if and only if $([h'_1]_F, ..., [h'_{\delta(Z)}]_F) \in Z^{\mathcal{A}_I/F}\textsubscript{-}$.
\end{enumerate}
\end{proposition}

\begin{prova}
We begin by 1.

\begin{enumerate}
\item $(\Rightarrow)$ Suppose $([h_1]_F, ..., [h_{\delta(Z)}]_F) \in Z^{\mathcal{A}_I/F}\textsubscript{+}$. Thus we have that $\{i \in I : (h_1(i), ..., h_{\delta(Z)}(i)) \in Z^{\mathcal{A}_i}\textsubscript{+}\} \in F$, and since $h_1 \sim_F h'_1, ..., h_{\delta(Z)} \sim_F h'_{\delta(Z)}$ we also have $\{i \in I: h_1(i) = h'_1(i)\}, ..., \{i \in I: h_{\delta(Z)}(i) = h'_{\delta(Z)}(i)\} \in F$, so that $\{i \in I: h_1(i) = h'_1(i)\} \cap ... \cap \{i \in I: h_{\delta(Z)}(i) = h'_{\delta(Z)}(i)\} \in F$. But then $\{i \in I: h_1(i) = h'_1(i)\} \cap ... \cap \{i \in I: h_{\delta(Z)}(i) = h'_{\delta(Z)}(i)\} \cap \{i \in I : (h_1(i), ..., h_{\delta(Z)}(i)) \in Z^{\mathcal{A}_i}\textsubscript{+}\} \in F$. Now if $k \in \{i \in I: h_1(i) = h'_{1}(i)\} \cap ... \cap \{i \in I: h_{\delta(Z)}(i) = h'_{\delta(Z)}(i)\} \cap \{i \in I : (h_1(i), ..., h_{\delta(Z)}(i)) \in Z^{\mathcal{A}_i}\textsubscript{+}\}$ then clearly $k \in \{i \in I : (h'_1(i), ..., h'_{\delta(Z)}(i)) \in Z^{\mathcal{A}_i}\textsubscript{+}\}$, whence $\{i \in I: h_1(i) = h'_1(i)\} \cap ... \cap \{i \in I: h_{\delta(Z)}(i) = h'_{\delta(Z)}(i)\} \cap \{i \in I : (h_1(i), ..., h_{\delta(Z)}(i)) \in Z^{\mathcal{A}_i}\textsubscript{+}\} \subseteq \{i \in I : (h'_1(i), ..., h'_{\delta(Z)}(i)) \in Z^{\mathcal{A}_i}\textsubscript{+}\}$. So $\{i \in I : (h'_1(i), ..., h'_{\delta(Z)}(i)) \in Z^{\mathcal{A}_i}\textsubscript{+}\} \in F$ and therefore $([h'_1]_F, ..., [h'_{\delta(Z)}]_F) \in Z^{\mathcal{A}_I/F}\textsubscript{+}$.

\medskip

$(\Leftarrow)$ The argument is similar.

\item The proof of item 2 is analogous to that of item 1.

\item $(\Rightarrow)$ (Contraposition) Assume that $([h'_1]_F, ..., [h'_{\delta(Z)}]_F) \not\in Z^{\mathcal{A}_I/F}\textsubscript{-}$. Thus either $([h'_1]_F, ..., [h'_{\delta(Z)}]_F) \in Z^{\mathcal{A}_I/F}\textsubscript{+}$ or $([h'_1]_F, ..., [h'_{\delta(Z)}]_F) \in Z^{\mathcal{A}_I/F}\textsubscript{0}$. So, by itens 1 and 2 it follows that $([h_1]_F, ..., [h_{\delta(Z)}]_F) \in Z^{\mathcal{A}_I/F}\textsubscript{+}$ or $([h_1]_F, ..., [h_{\delta(Z)}]_F) \in Z^{\mathcal{A}_I/F}\textsubscript{0}$ and hence $([h_1]_F, ..., [h_{\delta(Z)}]_F) \not\in Z^{\mathcal{A}_I/F}\textsubscript{-}$.

\medskip

$(\Leftarrow)$ The argument is similar.
\end{enumerate}

Therefore, the result holds.
\end{prova}

\begin{affirmation}
Consider the family \{$\mathcal{A}_i: i \in I$\}, and let the family \{$\mathcal{B}_i: i \in I$\} be such that $\mathcal{B}_i$ is $\mathcal{A}_i$-normal for each $i \in I$. Now let $Z \in \mathcal{R}$ and $h_1, ..., h_{\delta(Z)} \in A_I$ be any elements. Then the following assertions are true:

\begin{enumerate}
\item $\{i \in I : (h_1(i), ..., h_{\delta(Z)}(i)) \in Z^{\mathcal{A}_i}\textsubscript{+}\} \subseteq \{i \in I : (h_1(i), ..., h_{\delta(Z)}(i)) \in Z^{\mathcal{B}_i}\textsubscript{+}\}$;

\item $\{i \in I : (h_1(i), ..., h_{\delta(Z)}(i)) \in Z^{\mathcal{B}_i}\textsubscript{+}\} \subseteq \{i \in I : (h_1(i), ..., h_{\delta(Z)}(i)) \in Z^{\mathcal{A}_i}\textsubscript{+}\} \ \cup \ \{i \in I : (h_1(i), ..., h_{\delta(Z)}(i)) \in Z^{\mathcal{A}_i}\textsubscript{0}\}$;

\item $\{i \in I : (h_1(i), ..., h_{\delta(Z)}(i)) \in Z^{\mathcal{A}_i}\textsubscript{+}\} = \{i \in I : (h_1(i), ..., h_{\delta(Z)}(i)) \in Z^{{\mathcal{A}_i}_-}\textsubscript{+}\}$;

\item $\{i \in I : (h_1(i), ..., h_{\delta(Z)}(i)) \in Z^{{\mathcal{A}_i}_+}\textsubscript{+}\} = \{i \in I : (h_1(i), ..., h_{\delta(Z)}(i)) \in Z^{\mathcal{A}_i}\textsubscript{+}\} \cup \{i \in I : (h_1(i), ..., h_{\delta(Z)}(i)) \in Z^{\mathcal{A}_i}\textsubscript{0}\}$.
\end{enumerate}
\end{affirmation}

\begin{proposition}
Consider the family \{$\mathcal{A}_i: i \in I$\} and let \textit{F} be a filter over \textit{I}. Now let the family \{$\mathcal{B}_i: i \in I$\} be such that $\mathcal{B}_i$ is $\mathcal{A}_i$-normal for each $i \in I$. Then the following assertions are true:

\begin{enumerate}
\item $A_I/F = B_I/F$;

\item For $Z \in \mathcal{R}$, we have that $Z^{\mathcal{A}_I/F}\textsubscript{+} \subseteq Z^{\mathcal{B}_I/F}\textsubscript{+}$.
\end{enumerate}

\end{proposition}

\begin{prova}
We begin by 1.

\begin{enumerate}
\item Using the fact that $A_I = B_I$ by Proposition 6, we have:
\begin{eqnarray*}
[h]_{F} \in A_I/F &\Leftrightarrow& h \in A_I \hspace{2 cm} \text{(for each} \ [h]_F \in A_I/F)\\
&\Leftrightarrow& h \in B_I\\
&\Leftrightarrow& [h]_F \in B_I/F
\end{eqnarray*}

Hence, $A_I/F = B_I/F$.

\item For $Z \in \mathcal{R}$ and $[h_1]_F, ..., [h_{\delta(Z)}]_F \in A_I/F$, using item 1 of Affirmation 2, we have:
\begin{eqnarray*}
([h_1]_F, ..., [h_{\delta(Z)}]_F) \in Z^{\mathcal{A}_I/F}\textsubscript{+} &\Rightarrow& \{i \in I : (h_1(i), ..., h_{\delta(Z)}(i)) \in Z^{\mathcal{A}_i}\textsubscript{+}\} \in F\\
&\Rightarrow& \{i \in I : (h_1(i), ..., h_{\delta(Z)}(i)) \in Z^{\mathcal{B}_i}\textsubscript{+}\} \in F\\
&\Rightarrow& ([h_1]_F, ..., [h_{\delta(Z)}]_F) \in Z^{\mathcal{B}_I/F}\textsubscript{+}
\end{eqnarray*}

Hence, $Z^{\mathcal{A}_I/F}\textsubscript{+} \subseteq Z^{\mathcal{B}_I/F}\textsubscript{+}$.
\end{enumerate}

Therefore, the result holds.
\end{prova}

\begin{proposition}
Consider the family \{$\mathcal{A}_i: i \in I$\} and let \textit{F} be a filter over \textit{I}. Then, ${\mathcal{A}_I/F}_-$ is the reduced product of the family \{${\mathcal{A}_i}_- : i \in I$\} over \textit{F}.
\end{proposition}

\begin{prova}
Let $\mathcal{D} = (D, (Z^\mathcal{D})_{Z \in L})$ be the reduced product of \{${\mathcal{A}_i}_- : i \in I$\} over \textit{F} and let us show that ${\mathcal{A}_I/F}_- = \mathcal{D}$.

\begin{enumerate} [(a)]

\item We already have that $A_I/F = {A_I/F}_-$. Now by Proposition 9 we also have that $A_I/F = D$. Hence, ${A_I/F}_- = D$.

\item For $Z \in \mathcal{R}$ and $[h_1]_F, ..., [h_{\delta(Z)}]_F \in {A_I/F}_-$, using item 3 of Affirmation 4, we have:
\begin{eqnarray*}
([h_1]_F, ..., [h_{\delta(Z)}]_F) \in Z^{{\mathcal{A}_I/F}_-}\textsubscript{+} &\Leftrightarrow& ([h_1]_F, ..., [h_{\delta(Z)}]_F) \in Z^{\mathcal{A}_I/F}\textsubscript{+}\\
&\Leftrightarrow& \{i \in I : (h_1(i), ..., h_{\delta(Z)}(i)) \in Z^{\mathcal{A}_i}\textsubscript{+}\} \in F\\
&\Leftrightarrow& \{i \in I : (h_1(i), ..., h_{\delta(Z)}(i)) \in Z^{{\mathcal{A}_i}_-}\textsubscript{+}\} \in F\\
&\Leftrightarrow& ([h_1]_F, ..., [h_{\delta(Z)}]_F) \in Z^\mathcal{D}\textsubscript{+}
\end{eqnarray*}

Hence, $Z^{{\mathcal{A}_I/F}_-}\textsubscript{+} = Z^\mathcal{D}\textsubscript{+}$.

\end{enumerate}

Therefore, ${\mathcal{A}_I/F}_- = \mathcal{D}$.
\end{prova}

\bigskip

As we saw in Propositions 6 and 7, ${\mathcal{A}_I}_+$ is the direct product of the family \{${\mathcal{A}_i}_+ : i \in I$\}, and if a family \{$\mathcal{B}_i: i \in I$\} is such that $\mathcal{B}_i$ is $\mathcal{A}_i$-normal for each $i \in I$, the direct product of \{$\mathcal{B}_i: i \in I$\} is $\mathcal{A}_I$-normal. The next example shows, however, that we do not have similar results for reduced products. More specifically, if \textit{F} is a filter over \textit{I}, then $\mathcal{B}_I/F$ might not be $\mathcal{A}_I/F$-normal, and ${\mathcal{A}_I/F}_+$ might not be the reduced product of \{${\mathcal{A}_i}_+ : i \in I$\} over \textit{F}.

Assume that \{$\mathcal{A}_i: i \in I$\} is such that:

\begin{itemize}
\item $L = \mathcal{R} = \{R\}$;
\item $\delta(R) = 1$;
\item $I = \{x, y\}$;
\item $\mathcal{A}_x = (A_x, R^{\mathcal{A}_x})$ and $\mathcal{A}_y = (A_y, R^{\mathcal{A}_y})$;
\item $A_x = \{a_1\}$ and $A_y = \{a_2\}$;
\item $R^{\mathcal{A}_x}\textsubscript{+} = \{a_1\}, R^{\mathcal{A}_y}\textsubscript{0} = \{a_2\}$ and $R^{\mathcal{A}_x}\textsubscript{-} = R^{\mathcal{A}_x}\textsubscript{0} = R^{\mathcal{A}_y}\textsubscript{+} = R^{\mathcal{A}_y}\textsubscript{-} = \emptyset$.
\end{itemize}

\noindent Thus $\mathcal{A}_I$ is such that:

\begin{itemize}
\item $A_I = \{h\}$, where:

\begin{center}
$h : I \to A_x \cup A_y$

$x \mapsto {h}(x) = a_1$

$y \mapsto {h}(y) = a_2$;
\end{center}

\item $R^{\mathcal{A}_I}\textsubscript{0} = \{h\}$ and $R^{\mathcal{A}_I}\textsubscript{+} = R^{\mathcal{A}_I}\textsubscript{-} = \emptyset$.
\end{itemize}

\noindent Now suppose that $F = \{\{x, y\}\}$, i.e. \textit{F} is the trivial filter $\{I\}$. Then $\sim_F \ = {A_I}^2$, and as to $\mathcal{A}_I/F$ we have: 

\begin{itemize}
\item $A_I/F = \{[h]_F\}$;

\item $R^{\mathcal{A}_I/F}\textsubscript{-} = \{[h]_F\}$ and $R^{\mathcal{A}_I/F}\textsubscript{+} = R^{\mathcal{A}_I/F}\textsubscript{0} = \emptyset$.
\end{itemize}

\noindent Finally, \{${\mathcal{A}_i}_+ : i \in I$\} is such that:

\begin{itemize}
\item ${\mathcal{A}_x}_+ = ({A_x}_+, R^{{\mathcal{A}_x}_+})$ and ${\mathcal{A}_y}_+ = ({A_y}_+, R^{{\mathcal{A}_y}_+})$;
\item ${A_x}_+ = \{a_1\}$ and ${A_y}_+ = \{a_2\}$;
\item $R^{{\mathcal{A}_x}_+}\textsubscript{+} = \{a_1\}, R^{{\mathcal{A}_y}_+}\textsubscript{+} = \{a_2\}$ and $R^{{\mathcal{A}_x}_+}\textsubscript{-} = R^{{\mathcal{A}_x}_+}\textsubscript{0} = R^{{\mathcal{A}_y}_+}\textsubscript{-} = R^{{\mathcal{A}_y}_+}\textsubscript{0} = \emptyset$.
\end{itemize}

\noindent If we let $\mathcal{B} = (B, R^\mathcal{B})$ be the reduced product of \{${\mathcal{A}_i}_+ : i \in I$\} over \textit{F}, then:

\begin{itemize}
\item $B = \{[h]_F\} = A_I/F$;

\item $R^\mathcal{B}\textsubscript{+} = \{[h]_F\}$ and $R^\mathcal{B}\textsubscript{-} = R^\mathcal{B}\textsubscript{0} = \emptyset$.
\end{itemize}

As can be seen, $\mathcal{A}_I/F$ is total and $\mathcal{A}_I/F \neq \mathcal{B}$, so $\mathcal{B}$ is not $\mathcal{A}_I/F$-normal (Proposition 2).
Further, since $\mathcal{A}_I/F$ is total, $\mathcal{A}_I/F = {\mathcal{A}_I/F}_+$ and hence ${\mathcal{A}_I/F}_+ \neq \mathcal{B}$. The fact that $\mathcal{B}$ is not $\mathcal{A}_I/F$-normal and ${\mathcal{A}_I/F}_+ \neq \mathcal{B}$ holds only because if a set $X \cup Y$ belongs to a filter, it does not necessarily follow that either \textit{X} or \textit{Y} also belong. More precisely, if $\{x, y\} = \{i \in I: h(i) \in R^{\mathcal{A}_x}\textsubscript{+}\} \cup \{i \in I: h(i) \in R^{\mathcal{A}_y}\textsubscript{0}\} \in F$, it does not follow that $\{x\} = \{i \in I: h(i) \in R^{\mathcal{A}_x}\textsubscript{+}\} \in F$ or $\{y\} = \{i \in I: h(i) \in R^{\mathcal{A}_y}\textsubscript{0}\} \in F$. Otherwise we would have that $[h]_F \in R^{\mathcal{A}_I/F}\textsubscript{+}$ or $[h]_F \in R^{\mathcal{A}_I/F}\textsubscript{0}$, and $[h]_F \in R^{{\mathcal{A}_I/F}_+}\textsubscript{+}$, that is, we would have ${\mathcal{A}_I/F}_+ = \mathcal{B}$.

\subsection{Ultraproducts} 

We begin this third part of section 3 with the notion of ultrafilter, and we shall conclude it with the compactness theorem.

\begin{definition}
Let \textit{U} be a filter over \textit{I}. We say that \textit{U} is an \textit{ultrafilter over I} if for every $X \subseteq I$,

\begin{center}
$X \in U$ if and only if $I - X \not\in U$.
\end{center}
\end{definition}

\begin{affirmation}
Let \textit{U} be an ultrafilter over \textit{I}. For every $X, Y \subseteq I$,

\begin{center}
if $X \cup Y \in U$ then $X \in U$ or $Y \in U$.
\end{center}
\end{affirmation}

\begin{definition}
A collection \textit{E} of subsets of \textit{I} is said to have the \textit{finite intersection property} if and only if every intersection of finitely many elements of \textit{E} is non-empty.
\end{definition}

\begin{affirmation}
[Ultrafilter theorem] Every set with the finite intersection property can be extended to an ultrafilter.
\end{affirmation}

\begin{definition}
Let \textit{U} be a filter over \textit{I}, and consider the reduced product $\mathcal{A}_I/U$. We say that $\mathcal{A}_I/U$ is an \textit{ultraproduct}, if \textit{U} is an ultrafilter over \textit{I}. If $\mathcal{A}_I/U$ is the reduced power of an $\mathcal{L}$-partial structure $\mathcal{A}$ over \textit{U}, and \textit{U} is an ultrafilter, then $\mathcal{A}_I/U$ is said an \textit{ultrapower} of $\mathcal{A}$.
\end{definition}

As we have seen, if \textit{F} is a filter over \textit{I} then ${\mathcal{A}_I/F}_+$ might not be the reduced product of the family \{${\mathcal{A}_i}_+ : i \in I$\} over \textit{F}, and if \{$\mathcal{B}_i: i \in I$\} is a family such that $\mathcal{B}_i$ is $\mathcal{A}_i$-normal for each $i \in I$, then $\mathcal{B}_I/F$ might not be $\mathcal{A}_I/F$-normal. Now we will se that if \textit{F} is an ultrafilter those two results do not follow and this has to do with Affirmation 5 (which does not hold for filters in general, but only for ultrafilters).

\begin{proposition}
Consider the family \{$\mathcal{A}_i: i \in I$\} and let \textit{U} be an ultrafilter over \textit{I}. Now let the family \{$\mathcal{B}_i: i \in I$\} be such that $\mathcal{B}_i$ is $\mathcal{A}_i$-normal for each $i \in I$. Then $\mathcal{B}_I/U$ is ${\mathcal{A}_I/U}$-normal.
\end{proposition}

\begin{prova}
We already have that $\mathcal{B}_I/F$ is total, thus by Proposition 9 it only remains to verify that for $Z \in \mathcal{R}$, $Z^{\mathcal{A}_I/U}\textsubscript{-} \subseteq Z^{\mathcal{B}_I/U}\textsubscript{-}$. So, assuming that $[h_1]_U, ..., [h_{\delta(Z)}]_U \in A_I/U$ are any elements, and using item 2 of Affirmation 4 as well as Affirmation 5, we have:
\begin{eqnarray*}
\hspace{-3 cm} ([h_1]_U, ..., [h_{\delta(Z)}]_U) \in Z^{\mathcal{A}_I/U}\textsubscript{-} &\Rightarrow& \{i \in I : (h_1(i), ..., h_{\delta(Z)}(i)) \in Z^{\mathcal{A}_i}\textsubscript{+}\} \not\in U 
\ \text{and} \ \{i \in I : (h_1(i), ..., h_{\delta(Z)}(i)) \in Z^{\mathcal{A}_i}\textsubscript{0}\} \not\in U\\ 
&\Rightarrow& \{i \in I : (h_1(i), ..., h_{\delta(Z)}(i)) \in Z^{\mathcal{A}_i}\textsubscript{+}\} \cup \{i \in I : (h_1(i), ..., h_{\delta(Z)}(i)) \in Z^{\mathcal{A}_i}\textsubscript{0}\} \not\in U\\
&\Rightarrow& \{i \in I : (h_1(i), ..., h_{\delta(Z)}(i)) \in Z^{\mathcal{B}_i}\textsubscript{+}\} \not\in U\\
&\Rightarrow& ([h_1]_U, ..., [h_{\delta(Z)}]_U) \in Z^{\mathcal{B}_I/U}\textsubscript{-}
\end{eqnarray*}

Therefore, $\mathcal{B}_I/U$ is ${\mathcal{A}_I/U}$-normal.
\end{prova}

\begin{proposition}
${\mathcal{A}_I/U}_+$ is the reduced product of the family \{${\mathcal{A}_i}_+ : i \in I$\} over \textit{U}.
\end{proposition}

\begin{prova}
Suppose that $\mathcal{B} = (B, (Z^\mathcal{B})_{Z \in L})$ is the reduced product of \{${\mathcal{A}_i}_+ : i \in I$\} over \textit{U} and let us check that ${\mathcal{A}_I/U}_+ = \mathcal{B}$.

\begin{enumerate} [(a)]

\item The argument is similar to that of item (a) of Proposition 10.

\item For $Z \in \mathcal{R}$ and $[h_1]_U, ..., [h_{\delta(Z)}]_U \in {A_I/U}_+$, using item 4 of Affirmation 4 and Affirmation 5, we have:
\begin{eqnarray*}
\hspace{-4 cm} ([h_1]_U, ..., [h_{\delta(Z)}]_U) \in Z^{{\mathcal{A}_I/U}_+}\textsubscript{+} &\Leftrightarrow& ([h_1]_U, ..., [h_{\delta(Z)}]_U) \in Z^{\mathcal{A}_I/U}\textsubscript{+} \cup Z^{\mathcal{A}_I/U}\textsubscript{0}\\
&\Leftrightarrow& ([h_1]_U, ..., [h_{\delta(Z)}]_U) \in Z^{\mathcal{A}_I/U}\textsubscript{+} \ \text{or} \ ([h_1]_U, ..., [h_{\delta(Z)}]_U) \in Z^{\mathcal{A}_I/U}\textsubscript{0}\\
&\Leftrightarrow& \{i \in I : (h_1(i), ..., h_{\delta(Z)}(i)) \in Z^{\mathcal{A}_i}\textsubscript{+}\} \in U \ \text{or} \ \{i \in I : (h_1(i), ..., h_{\delta(Z)}(i)) \in
 Z^{\mathcal{A}_i}\textsubscript{0}\} \in U\\
&\Leftrightarrow& \{i \in I : (h_1(i), ..., h_{\delta(Z)}(i)) \in Z^{\mathcal{A}_i}\textsubscript{+}\} \cup \{i \in I : (h_1(i), ..., h_{\delta(Z)}(i)) \in Z^{\mathcal{A}_i}\textsubscript{0}\} \in U\\
&\Leftrightarrow& \{i \in I : (h_1(i), ..., h_{\delta(Z)}(i)) \in Z^{{\mathcal{A}_i}_+}\textsubscript{+}\} \in U\\
&\Leftrightarrow& ([h_1]_U, ..., [h_{\delta(Z)}]_U) \in Z^\mathcal{B}\textsubscript{+}
\end{eqnarray*}

Hence, $Z^{{\mathcal{A}_I/U}_+}\textsubscript{+} = Z^\mathcal{B}\textsubscript{+}$.
\end{enumerate}

Therefore, ${\mathcal{A}_I/U}_+ = \mathcal{B}$.
\end{prova}

\begin{lemma}
Let $\varphi$ be an $\mathcal{L}$-sentence. The following assertions are true:

\begin{enumerate}
\item If \{$\mathcal{B}_i: i \in I$\} is a family such that $\mathcal{B}_i$ is $\mathcal{A}_i$-normal for each $i \in I$, then $\{i \in I: \mathcal{B}_i \models \varphi\} \subseteq \{i \in I: \mathcal{A}_i \mid\models \varphi\}$;

\item There exists a family \{$\mathcal{B}_i: i \in I$\} such that $\mathcal{B}_i$ is $\mathcal{A}_i$-normal for each $i \in I$, and $\{i \in I: \mathcal{A}_i \mid\models \varphi\} = \{i \in I: \mathcal{B}_i \models \varphi\}$.
\end{enumerate}

\end{lemma}

\begin{prova}
Assume that $\varphi$ is an $\mathcal{L}$-sentence. 

\begin{enumerate}
\item Suppose that the family \{$\mathcal{B}_i: i \in I$\} is such that $\mathcal{B}_i$ is $\mathcal{A}_i$-normal for each $i \in I$. Then:
\begin{eqnarray*}
k \in \{i \in I: \mathcal{B}_i \models \varphi\} &\Rightarrow& \mathcal{B}_k \models \varphi \hspace{3 cm} \text{(for each $k \in \{i \in I: \mathcal{B}_i \models \varphi\}$)}\\
&\Rightarrow& \mathcal{A}_k \mid\models \varphi\\
&\Rightarrow& k \in \{i \in I: \mathcal{A}_i \mid\models \varphi\}
\end{eqnarray*}

\item Define the family \{$\mathcal{B}_i: i \in I$\} as follows:

\begin{enumerate} [i.]
\item If $k \in \{i \in I: \mathcal{A}_i \mid\models \varphi\}$, then $\mathcal{B}_k$ is $\mathcal{A}_k$-normal and $\mathcal{B}_k \models \varphi$;

\item If $k \not\in \{i \in I: \mathcal{A}_i \mid\models \varphi\}$, then $\mathcal{B}_k$ is any $\mathcal{A}_k$-normal structure.
\end{enumerate}

It is immediate that $\{i \in I: \mathcal{A}_i \mid\models \varphi\} \subseteq \{i \in I: \mathcal{B}_i \models \varphi\}$, and $\{i \in I: \mathcal{B}_i \models \varphi\} \subseteq \{i \in I: \mathcal{A}_i \mid\models \varphi\}$ by item 1. Hence, $\{i \in I: \mathcal{A}_i \mid\models \varphi\} = \{i \in I: \mathcal{B}_i \models \varphi\}$.
\end{enumerate}

Therefore, the result holds. 
\end{prova}

\begin{affirmation} 
[Ło$\check{s}$'s theorem] Let \textit{U} be an ultrafilter over \textit{I}. Let \{$\mathcal{B}_i: i \in I$\} be a family of total $\mathcal{L}$-structures and consider the ultraproduct $\mathcal{B}_I/U$. For each $\mathcal{L}$-sentence $\varphi$,

\begin{center}
$\mathcal{B}_I/U \models \varphi$ if and only if $\{i \in I: \mathcal{B}_i \models \varphi\} \in U$.
\end{center}
\end{affirmation}

Next we prove part of Ło$\check{s}$'s theorem with respect to partial structures and quasi-truth. This part is the only one we will need to prove the compactness theorem.

\begin{proposition}
Let \textit{U} be an ultrafilter over \textit{I} and consider the ultraproduct $\mathcal{A}_I/U$. For each $\mathcal{L}$-sentence $\varphi$,

\begin{center}
if $\{i \in I: \mathcal{A}_i \mid\models \varphi\} \in U$ then $\mathcal{A}_I/U \mid\models \varphi$.
\end{center}
\end{proposition}

\begin{prova}
Let $\varphi$ be an $\mathcal{L}$-sentence and assume that $\{i \in I: \mathcal{A}_i \mid\models \varphi\} \in U$. By item 2 of Lemma 3, there exists a family \{$\mathcal{B}_i: i \in I$\} such that $\mathcal{B}_i$ is $\mathcal{A}_i$-normal for each $i \in I$, and $\{i \in I: \mathcal{A}_i \mid\models \varphi\} = \{i \in I: \mathcal{B}_i \models \varphi\}$. Hence $\{i \in I: \mathcal{B}_i \models \varphi\} \in U$, and so $\mathcal{B}_I/U \models \varphi$ by Affirmation 7. Now given that $\mathcal{B}_I/U$ is $\mathcal{A}_I/U$-normal by Proposition 11, we have $\mathcal{A}_I/U \mid\models \varphi$.  
\end{prova}

\begin{proposition}
[compactness theorem] A set of $\mathcal{L}$-sentences $\Gamma$ has a partial model if and only if every finite subset of $\Gamma$ has a partial model.
\end{proposition}

\begin{prova} 
$(\Rightarrow)$ Straightforward.

\medskip

$(\Leftarrow)$ Let $\Gamma$ be a set of $\mathcal{L}$-sentences, \textit{I} be the set of all finite subsets of $\Gamma$, and suppose that for every $i \in I$, there exists a partial model $\mathcal{A}_i$ of \textit{i}. Now for every $\gamma \in \Gamma$, let $\gamma^*$ be the set of all $i \in I$ such that $\gamma \in i$, and 

\begin{center}
$E = \{\gamma^* : \gamma \in \Gamma\}$.
\end{center}

Then it is easy to check that \textit{E} has the finite intersection property, whence it can be extended to an ultrafilter over \textit{I} by Affirmation 6 (ultrafilter theorem). So let \textit{U} be an ultrafilter over \textit{I} such that $E \subseteq U$. Thus for each $\gamma \in \Gamma, \gamma^* \in U$. Moreover, $\gamma^* \subseteq \{i \in I:\mathcal{A}_i \mid\models \gamma\}$, since if $i \in \gamma^*$ then $\gamma \in i$ and hence $\mathcal{A}_i \mid\models \gamma$. So, for each $\gamma \in \Gamma$, it follows that $\{i \in I: \mathcal{A}_i \mid\models \gamma\} \in U$. But then, for each $\gamma \in \Gamma$, we have that $\mathcal{A}_I/U \mid\models \gamma$ by Proposition 13. Therefore, $\mathcal{A}_I/U \mid\models \Gamma$.
\end{prova}

\section{Concluding remarks}

The content of section 3 is relevant for two main reasons: because of the compactness theorem, firstly, but also because it shows that an important \textit{method of constructing structures} can be applied to partial structures, and used to prove the result mentioned. More specifically, we saw that one can define the notion of ultraproducts (just like the notions of direct product and reduced product) for partial structures, and use it to demonstrate the compactness theorem. Accordingly, partial model theory preserves more than this result, despite the fact that its logic is not the same as that of traditional (or standard) model theory, i.e. classical logic.

In addition to analyze whether other results involving the notions of reduced product, direct product and ultraproduct are preserved in partial model theory, the next step in the development of this theory is to analyze whether other methods of constructing structures could be applied to partial structures, like the method of \textit{Elementary chains} (cf. \cite{Chang1990} p. v), and whether other fundamental results of traditional model theory are preserved too, like the \textit{Löwenhein-Skolem theorem}. This is something we intend to undertake in a series of further papers.

\section*{Acknowledgments} This work was supported by The São Paulo Research Foundation (FAPESP) under Grant Agreement No 2022/01401-8.

\end{document}